\documentclass[11pt]{article}
\usepackage{latexsym,amssymb,amsmath}
\input xy
\xyoption{all}
\def\fa{\mathfrak a}

\def\cwedge{\wedge\cdots\wedge}

\def\trace{{\rm trace}}

\def\cS{{\mathcal S}}

\def\cO{{\mathcal O}}

\def\CC{\mathbb C}
\def\RR{\mathbb R}
\def\HH{\mathbb H}
\def\AA{{\mathbb A}}

\def\OO{\mathbb O}

\def\SS{\mathbb S}

\def\11{\mathbf 1}
\def\PP{\mathbb P}

\def\FF{\mathbb F}

\def\e1{\varepsilon_1}
\def\e2{\varepsilon_2}
\def\e3{\varepsilon_3}

\def\P2{{\PP}^2}

\def\JA{{\cal J}_3(\AA)}

\def\00{\underline{0}}
\def\J0{{\cal J}_3(\underline{0})}

\def\PJ0{\PP({\cal J}_3(\underline{0}))}

\def\fsl{{\mathfrak {sl}}}
\def\fgl{{\mathfrak {gl}}}

\def\fc{{\mathfrak c}}
\def\fa{{\mathfrak a}}

\def\ft{{\mathfrak t}}

\def\a{\alpha}

\def\om{\omega}

\def\b{\beta}
\def\g{\gamma}
\def\s{\sigma}

\def\d{\delta}

\def\e{\varepsilon}
\def\ot{{\mathord{\,\otimes }\,}}
\def\op{{\mathord{\,\oplus }\,}}

\def\lra{{\mathord{\;\longrightarrow\;}}}
\def\ra{{\mathord{\;\rightarrow\;}}}

\def\we{{\mathord{{\scriptstyle \wedge}}}}

\def\AP2{{\AA\PP}^2}
\def\RP2{{\RR\PP}^2}
\def\CP2{{\CC\PP}^2}
\def\HP2{{\HH\PP}^2}
\def\OP2{{\OO\PP}^2}

\def\tr{{\rm trace}\;}
\def\dim{{\rm dim}\;}

\newcommand\proof{{\noindent {\em Proof}.}\hspace{2mm}}
\newcommand\qed{{\hfill\hfill $\Box$}}

\newtheorem{theo}{Theorem}
\newtheorem{coro}[theo]{Corollary}
\newtheorem{lemm}[theo]{Lemma}
\newtheorem{prop}[theo]{Proposition}

\begin{document}

\title{Varieties of reductions for $\fgl_n$}
\author{A. Iliev, L. Manivel}
\maketitle

\begin{abstract} We study the varieties of reductions associated to the
variety of rank one matrices in $\fgl_n$. In particular, we prove that 
for $n=4$ we get a $12$-dimensional Fano variety of Picard number one
 and index $3$, with canonical singularities.
\end{abstract}

\section{Introduction}

This paper is a sequel to \cite{imcrelle} and the companion 
paper \cite{imcomp},
where we studied a family of smooth Fano varieties with many remarkable
properties. These varieties were constructed as compactifications of 
what we called {\it reductions} for the four Severi varieties. Recall 
that the  Severi varieties can be defined as the 
projective planes over the 
four (complexified) normed algebras $\AA=\RR,\CC,\HH,\OO$ -- the reals, 
the complexes, the quaternions, and the octonions. More precisely,
consider the Jordan algebra $\JA$ of $\AA$-Hermitian matrices of order
$3$. The projectivization of the set of rank one matrices in $\PP\JA$ 
is the Severi variety $X_a$, a homogeneous variety of dimension $2a$, 
where $a=1,2,4,8$ denotes the dimension of $\AA$.
 
A {\it non singular reduction} is defined as a $3$-secant plane 
to $X_a$ passing through the identity matrix $I$. The projection 
$p$ from $I$
to the hyperplane $\PP\JA_0$ of traceless matrices sends the 
non-singular reductions to the family of $3$-secant lines 
to the projected Severi variety $\overline{X}_a$, and the 
variety of reductions that we studied in \cite{imcrelle} is the
compactification of that family in the Grassmannian of lines in 
$\PP\JA_0$. We proved that it is a smooth Fano manifold of dimension
$3a$, Picard number one, and index $a+1$.

In this paper we consider matrices of rank order than three, 
and the corresponding varieties of reductions. For $a=1$ they were
previously studied by Ranestad and Schreyer \cite{rs}, who proved 
that they are smooth up to rank $5$, while in rank $6$ the tangent
cone to a normal slice to the singular locus is, rather remarkably, 
a cone over the spinor variety $\SS_{10}$. Here we will focus on the 
case $a=2$, which has the interesting feature of being related to 
different, but not less classical problems than the study of Fano
varieties. Indeed, a non singular reduction for the variety of rank
one matrices $X_{2,n}=\PP^{n-1}\times {\check\PP}^{n-1}\subset\PP\fgl_n$,
is the commutative algebra of matrices that are diagonal with respect
to some basis of $\CC^n$ -- hence a direct connection with the
much studied problem of classifying commutative subalgebras of 
$\fgl_n$. Also, our variety of reductions $Red(n)$ appears as a 
natural compactification of the homogeneous space $PGL_n/N$, where 
$N$ denotes the normalizer of a maximal torus. 

For arbitrary $n$ a deep understanding of this compactification 
remains out of our reach : we only establish rather basic properties
and raise a number of questions. We mainly prove that $Red(n)$ 
is smooth in codimension one but always singular for $n\ge 4$. 
Moreover, the canonical divisor of the smooth locus is minus three
times the hyperplane divisor -- but we don't know if our varieties 
of reductions are normal in general. A tempting way to study 
$Red(n)$ is to consider its tautological fibration, which is 
birational to $\PP\fsl_n$. Quite interestingly, the induced 
rational map from this space to $Red(n)$ is closely related to 
the geometry of the set of non regular matrices. We only sketch
what should be the relevant plethystic transformations, and the 
connection with the Hilbert scheme of $n$ points in $\PP^{n-1}$. 

We can say a lot more when $n=4$. We prove that every abelian 
four dimensional subalgebra of $\fgl_4$ is in $Red(4)$, which is 
made of fourteen $PGL_4$-orbits. Three of these are closed, among 
which a projective three-space and its dual constitute the 
singular locus of $Red(4)$. We prove that the tangent cone to a normal
slice to each of these singular components is a cone over the 
Grassmannian $G(2,6)$ -- in particular, $Red(4)$ is normal. 
Blowing them up, we get a smooth variety in which a maximal torus 
of $PGL_4$ only has a finite number of fixed points. This allows 
us to compute the ranks of the Chow groups of $Red(4)$. We conclude 
that $Red(4)$ is a rational Fano variety of dimension $12$, 
Picard number one, index $3$, with canonical singularities.

Of course we expect that the variety of reductions defined for the
quaternions have similar properties, the geometry of the {\it Scorza
varieties} being quite insensitive to the underlying normed algebra
(see e.g. \cite{chaput}).

\section{Reductions for $\fgl_n$}

\subsection{Reductions and abelian algebras}
Let $Red(n)^0\subset G(n-1,\fsl_n)$ denote the space of Cartan
subalgebras of $\fsl_n$. Recall that $PGL_n$ acts transitively on Cartan
subalgebras, which are just the algebras of diagonal matrices 
with respect to some basis. Of course we may (and we will freely)
identify them with 
Cartan subalgebras  of $\fgl_n$, one way by adding the identity 
matrix, the other way by the natural projection $p : \fgl_n\ra\fsl_n$ 
from the identity matrix. From the point of view of reductions, 
a Cartan subalgebra of $\fgl_n$ is seen as a $n$-secant linear space to
the rank one variety
$X_n=X_{2,n}=\PP^{n-1}\times {\check\PP}^{n-1}\subset\PP\fgl_n$. 
Indeed, if such a linear space meets $X_n$ at $n$ distinct points
$e_1^*\ot e_1, \ldots, e_n^*\ot e_n$, and passes through $I$, we may
suppose that $I=e_1^*\ot e_1+ \cdots +e_n^*\ot e_n$, and then
automatically $e_1,\ldots ,e_n$ is a basis and $e_1^*,\ldots ,e_n^*$
is the dual basis. 

Once a Cartan subalgebra $\fa$ of $\fsl_n$ is fixed, we get
an isomorphism of $Red(n)^0$ with $PGL_n/N(\fa)$, where the normalizer
$N(\fa)$ is an extension of the maximal torus $A\subset PGL_n$ whose
Lie algebra is $\fa$, by the symmetric group $\cS_n$.
Let $Red(n)$ be the Zariski closure of $Red(n)^0$ in the Grassmannian
$G(n-1,\fsl_n)$. This compactification of $PGL_n/N(\fa)$ will be our
main object of interest. We call it the {\it variety of reductions}
for $\fsl_n$ (or $\fgl_n$).

First note that $Red(n)$ is a subvariety of the space $Ab(n)$ 
of abelian $(n-1)$-dimensional subalgebras of $\fsl_n$. 
This variety $Ab(n)$ has a simple set-theoretical description as the 
intersection of  $G(n-1,\fsl_n)\subset\PP\Lambda^{n-1}\fsl_n$ 
with the (projectivised) kernel of the natural map
$$\Theta :\Lambda^{n-1}\fsl_n\hookrightarrow \Lambda^2\fsl_n \ot
\Lambda^{n-3}\fsl_n\lra \fsl_n\ot \Lambda^{n-3}\fsl_n,$$
where the first arrow is the natural inclusion, and the second one is 
induced by the Lie bracket.  

Beware that this intersection is not transverse, and even not 
proper already for $n=3$, although $Red(3)$ turns out to be smooth. 
Moreover, $Ab(3)=Red(3)$, and we'll prove in the second part of this 
paper  that $Ab(4)=Red(4)$. An easy general result is:

\begin{prop}
The variety of reductions $Red(n)$ is an irreducible component of $Ab(n)$.
\end{prop}

\proof The generic element of a maximal torus in $\fsl_n$ is a semisimple 
endomorphism with distinct eigenvalues. Since having distinct eigenvalues 
is an open condition in $\fsl_n$, containing such an endomorphism is also 
an open condition in $Ab(n)$. But an abelian subalgebra of dimension 
$n-1$ in $\fsl_n$, which contains an endomorphism with distinct eigenvalues, 
must be the centralizer of this endomorphism -- hence a Cartan subalgebra.
This proves our claim. \qed

\medskip
In fact it is easy to show that $Ab(n)\neq Red(n)$ for large $n$. 
For example, suppose that $n=2m$ and let $L$ be any subspace of dimension 
$m$ in $\CC^n$. Let $\fa(L)$ denote the space of endomorphisms whose image is 
contained in $L$ and whose kernel contains $L$. Its dimension is $m^2$, and 
any $(n-1)$-dimensional subspace of  $\fa(L)$ is an abelian subalgebra of 
$\fsl_n$. Since a generic such subspace determines $L$ uniquely, we get 
a family of dimension $m^2+(2m-1)(m-1)^2$ in $A(n)$, which is strictly 
bigger than the dimension $n(n-1)$ of $Red(n)$ as soon as $m\ge 4$. 
A variant leads to the same conclusion for $n=2m-1$ and $m\ge 4$. 

\medskip\noindent {\it Question A}. 
Does $Ab(n)=Red(n)$ for $n=5$ or $6$?

\medskip\noindent {\it Remark}. 
 Suprunenko and Tyshkevich \cite{supru}
described explicitly the maximal nilpotent and
abelian subalgebras of $\fsl_5$ and $\fsl_6$. It turns out that there
is only a finite number of them up to conjugation, while there exists an
infinity in $\fsl_n$, $n\ge 7$. In principle this should allow to answer 
the question above. Indeed, by the Jordan decomposition, the semisimple
parts of the elements of an abelian subalgebra of $\fsl_n$ commute, so
that we can find a minimal decomposition of $\CC^n$ preserved by 
these, and basically, if this decomposition is not trivial, 
we are reduced to $\fsl_m$ with $m<n$. If the decomposition is trivial,
our subalgebra is nilpotent and we can use Suprunenko's results. 
 
\medskip 
The description of the other irreducible components of $Ab(n)$ is 
certainly an interesting problem. A basic question about the variety 
of reductions is:

\medskip\noindent {\it Question B}. How can we characterize the 
points of $Red(n)$ among the abelian subalgebras ? 

\medskip\noindent {\it Remark}. 
A necessary condition for an abelian algebra $\fa\in Ab(n)$ to belong to 
$Red(n)$, is that the commutative subalgebra of $\fgl_n$, generated by $\fa$
for the usual matrix product, has dimension at most $n$ (this was 
already pointed out by Gerstenhaber \cite{gerst}). 
But we don't know any example of 
an abelian subalgebra in $Ab(n)$ which does not fulfill this condition. 

\medskip Our hope is that $Red(n)$ should in general be a much 
nicer variety than $Ab(n)$ or its other irreducible components, when
they exist. For example, we observe that:

\begin{prop} The action of $PGL_n$ on $Ab(n)$ has finitely many orbits 
only for $n\le 5$. 
\end{prop}

\proof For $n\le 5$ this follows from the work of Suprunenko
and Tyshkevich \cite{supru}. 
Now suppose that $n\ge 6$, and that $n=2m$ is even. As above, let 
$L$ be an $m$-dimensional and consider the space of endomorphisms 
$\fa(L)$. Any $(n-1)$-dimensional subspace of  $\fa(L)$ is an 
abelian subalgebra of $\fsl_n$, and a  
 generic such subspace determines $L$. For $PGL_n$ to have 
a finite number of orbits in $Ab(n)$, the parabolic subgroup $P_L$ of $PGL_n$ 
stabilizing $L$ must have a finite number of orbits on the open subset 
$G(n-1,\fa(L))^0$ of  $(n-1)$-dimensional subspaces of  $\fa(L)$ whose generic
element has image $L$. But $P_L$ acts on the Grassmannian 
$G(n-1,\fa(L))$ only 
through its semisimple part $PGL_m\times PGL_m$, whose action is equivalent 
to its natural action on $G(n-1,M_m(\CC))$. The dimension of this 
Grassmannian
is strictly bigger than the dimension of  $PGL_m\times PGL_m$ as soon 
as $m\ge 3$, 
so there must be an infinity of orbits on any open subset. 

The case of odd $n$ is similar. \qed 

\medskip
In particular, the action of $PGL_n$ on $Red(n)$ has finitely many orbits 
for $n\le 5$. 

\medskip\noindent {\it Question C}. Does $Red(n)$ contain infinitely
many orbits of $PGL_n$ for $n\ge 6$ ?

\subsection{Special orbits}

A point in $Red(n)^0$ can be described as the centralizer of
a regular semisimple element of $\fsl_n$. If we drop the semisimplicity
hypothesis, we still get abelian $(n-1)$-dimensional subalgebras 
of $\fsl_n$ which we call {\it one-regular subalgebras}. Such 
subalgebras belong to $Red(n)$, as follows from the proof of our next 
result. 

\begin{prop}
The variety of reductions $Red(n)$ contains a unique codimension 
one orbit $\cO_{bound}$. A point in $\cO_{bound}$ is the centralizer 
of a regular matrix whose semisimple part has an eigenvalue of
multiplicity two.
\end{prop}

\proof 
Indeed, a point in this set $\cO_{bound}$ is defined by $n-1$ 
points in $\PP^{n-1}$, plus a plane containing one of the lines, 
all these spaces being in general position -- in particular, 
$PGL_n$ acts transitively on $\cO_{bound}$. Counting dimensions, 
we easily check that its codimension in $Red(n)$ equals one. 

Now let $\fa$ be a point of $Red(n)-Red(n)^0$, and consider a 
general point $x$ in $\fa$. By hypothesis, $x$ is not regular
semi-simple. Thus it belongs to the closure of the set of regular 
non-semisimple elements of $\fsl_n$. But this implies that $\fa$
belongs to the closure of the set of centralizers of such elements, 
thus to the closure of $\cO_{bound}$. In particular, if $\fa$
does not belong to $\cO_{bound}$, it must belong to a $PGL_n$-orbit
of smaller dimension. \qed

\medskip Extending this a little bit we can describe other 
orbits in $Red(n)$. Call an algebra $\fa\in Ab(n)$ 
{\it two-regular} if it can
be defined as the common centralizer of two of its elements. The 
irreducibility of the commuting variety \cite{mt} implies:

\begin{prop}
Any one or two-regular algebra in $Ab(n)$ does belong to $Red(n)$.
\end{prop}

On the other hand we can describe lots of closed orbits in $Ab(n)$.
If we choose a flag of subspaces of $\CC^n$, of the form
$$V_{i_1}\subset\cdots\subset V_{i_p}\subset V_{j_0}\subset V_{j_1}
\subset\cdots\subset V_{j_p},
$$
and if we consider the set of endomorphisms of $\CC^n$ mapping 
$V_{j_k}$ to $V_{i_k}$ for $k=1,\ldots,p$, and mapping $\CC^n$ to
$V_{j_0}$ and $V_{j_0}$ to zero, 
we get an abelian subalgebra of $\fsl_n$, which belongs to 
$Ab(n)$ when it has the correct dimension, that is, when
$$n-1=\sum (j_k-j_{k-1})(i_l-i_{l-1}).$$
When the flag varies, we get a closed $PGL_n$-orbit in $Ab(n)$,
but it is not clear to us whether it belongs to $Red(n)$ or not. 

A simple example is the case where our flag reduces to $V_{j_0}$,
which needs to be either a line or a hyperplane for the dimension
condition to be fulfilled. We thus get two
closed orbits $\cO_{min}'\simeq\PP^{n-1}$ and $\cO_{min}''\simeq {\check
\PP}^{n-1}$, which are dual projective spaces. 

\begin{prop}
The orbits $\cO_{min}'$ and $\cO_{min}''$ are contained in  $Red(n)$.
\end{prop}

\proof Consider the algebra of diagonal matrices with respect to a basis 
of the form $e_1,e_1+te_2,\ldots ,e_1+te_n$, and let $t$ tends to zero.
An easy computation shows that the limit point in $G(n-1,\fsl_n)$ 
belongs to $\cO_{min}'$, which is thus contained in $Red(n)$ -- hence
$\cO_{min}''$ as well, by duality. \qed

\subsection{Smoothness}

For $n\ge 4$, the variety of reductions  $Red(n)$ will be singular, 
but we expect the singular locus to be relatively small. Our main
general result in that direction is the following:

\begin{prop}
The codimension one orbit $\cO_{bound}$ is contained in the smooth 
locus of $Red(n)$. 

In particular $Red(n)$ is smooth in codimension one.
\end{prop}

\proof We choose a representative of $\cO_{bound}$ by fixing 
a basis $e_1,\ldots ,e_n$ of $\CC^n$ and letting 
$$\phi_1=e_2^*\ot e_1, \quad \phi_2=e_1^*\ot e_1+e_2^*\ot e_2,
\quad \phi_k=e_k^*\ot e_k \;\; {\rm for}\;\; 2<k\le n.$$
We check by an explicit computation that the Zariski tangent 
space to $Ab(n)$ at this point has dimension $n(n-1)$. Moreover,
a first order deformation in $Ab(n)$ (or $Red(n)$) is given, in
matrices, by

$$\psi_1=\begin{pmatrix} \mu & 1 &-\theta_{23} & \cdots & -\theta_{2n} \\
\nu & -\mu  & 0 & \cdots & 0 \\
0 & -\theta_{31}  & 0 & \cdots & 0 \\ \cdots & \cdots  & \cdots & \cdots & \cdots \\
0 & -\theta_{n1}  & 0 & \cdots & 0 \end{pmatrix},
\qquad
\psi_2=\begin{pmatrix} 1 & 0 & -\theta_{13} & \cdots 
& -\theta_{1n} \\ 0 & 1  & -\theta_{23} & \cdots 
& -\theta_{2n} \\ -\theta_{31} & -\theta_{32}  & 0 & \cdots & 0 \\
\cdots & \cdots  & \cdots & \cdots & \cdots \\
-\theta_{n1} & -\theta_{n2}  & 0 & \cdots & 0 \end{pmatrix},$$
and for $3\le k\le n$, 
$$\psi_k=\begin{pmatrix} 0 & 0 & \cdots &\cdots &  \theta_{1k} & \cdots & \cdots \\
0 & 0 & \cdots &\cdots &  \theta_{2k} &\cdots &  \cdots \\
\cdots & \cdots  & \cdots & \cdots & \theta_{3k} & \cdots & \cdots \\
\cdots & \cdots  & \cdots & \cdots & \cdots & \cdots & \cdots \\
\theta_{k1} & \theta_{k2} & \theta_{k3}  & \cdots & 1 & \cdots & 
-\theta_{kn} \\
\cdots & \cdots  & \cdots & \cdots & \cdots & \cdots & \cdots \\
\cdots & \cdots  & \cdots & \cdots & -\theta_{nk} & \cdots & \cdots \\
 \end{pmatrix}.$$
(The matrix $\psi_k$ has non zero coefficients only on the $k$-th line
and $k$-th column. Note the change of sign after the diagonal $1$.)
The number of free coefficients is
$2\binom{n-2}{2}+4(n-2)+2=n(n-1)$, as it should be. \qed

\medskip
We have a very simple geometric description of 
$\cO_{bound}$, which will be useful later.

\begin{prop}
The closure of $\cO_{bound}$ is a generically transverse quadric 
section of $Red(n)$. \end{prop}

\proof The Killing form induces a $PGL_n$-invariant quadric hypersurface 
in $\PP\Lambda^{n-1}\fsl_n$, given by 
$$Q(X_1\cwedge X_{n-1})=\det (\trace(X_iX_j))_{1\le i,j\le n-1}.$$
Note that this is the restriction of the quadric in
$\PP\Lambda^n\fgl_n$ given by (almost) the same formula,
the embedding being given by the wedge product with $I$.

Clearly this quadric does not contain $Red(n)^0$ but does contain 
its boundary. To check that the intersection is generically transverse
we compute $Q$ to first order on the matrices above. At first order, 
$\trace \psi_1^2=2\nu$, $\trace \psi_2^2=2$, $\trace \psi_k^2=1$
for $k>2$, and $\trace \psi_i\psi_j=0$ for $i\ne j$. Hence
$Q(\psi_1,\ldots ,\psi_n)=4\nu$, which proves our claim. \qed

\medskip
The tangent space to $Red(n)$ at a generic point $\fa$ is 
the image of the adjoint action 
$$\fsl_n\stackrel{ad}{\lra} Hom(\fa,\fsl_n/\fa)=T_{\fa}G(n-1,\fsl_n),$$
whose kernel is the normalizer of $\fa$, that is, $\fa$ itself
at the generic point. 

Note that this makes sense for any $\fa$ which is its own normalizer, 
in particular for any point of a regular orbit in $Red(n)$. We deduce
that the reduced tangent cone at such a point is linear, of the dimension
of $Red(n)$. This does not quite prove that we get a smooth point 
of $Red(n)$, but we can ask:

\medskip\noindent {\it Question D}. Is the set of one-regular
subalgebras contained in the smooth locus of $Red(n)$ ?

\subsection{The canonical sheaf}

For $\fa\in Red(n)^0$, we may identify $\fsl_n/\fa$
with the orthogonal $\fa^{\perp}$ of $\fa$ with respect to the 
Killing form, hence $\det T_{\fa}Red(n)$ with $\we^{top}\fa^{\perp}$,
the maximal wedge power.
Note that the the maximal torus $A$ in $PGL_n$ whose Lie algebra is 
$\fa$ acts trivially on this line. We thus get an action of the Weyl group
$N(A)/A\simeq\cS_n$, which is simply given by the sign representation. 
We deduce that the square $K_{Red(n)^0}^2$ of the canonical line bundle
of $Red(n)^0$, is trivial. Indeed, we can choose an orthonormal basis
of $\fa^{\perp}$ with respect to the Killing form,
 and consider the square of the 
corresponding volume form on $T_{\fa}Red(n)$. Since it is left invariant 
by the stabilizer of $\fa$ in $PGL_n$, we can translate it by $PGL_n$
to get a well-defined non vanishing section $\om$ of $K_{Red(n)^0}^2$. 

\smallskip
Let us compute the vanishing order 
of this section along the codimension one orbit 
$\cO_{bound}$. To do this we restrict to the following line in $Red(n)$,
which meets $\cO_{bound}$ transversely at $t=0$:

$$\fa(t)=\Bigg\{ \begin{pmatrix}
a_2 & a_1 & & & \\ ta_1 & a_2 & & & \\
 & & a_3 & & \\  & & & \cdots & \\
 & & & & a_n \end{pmatrix}, \quad a_1,\ldots ,a_n\in\CC \Bigg\}.
$$
What we first need is a first order deformation of $\fa(t)$ in 
$Red(n)$ for each $t$. We claim that such a deformation is provided by 
the following matrices:

$$\psi_1=\begin{pmatrix} \mu & 1 & -\theta_{23} & \cdots & -\theta_{2n} \\
t+\nu & -\mu  & -t\theta_{13} & \cdots & -t\theta_{1n} \\
-t\theta_{32} & -\theta_{31}   & 0 & \cdots & 0 \\ \cdots & \cdots
  & \cdots & \cdots & \cdots \\
-t\theta_{n2} & -\theta_{n1}  & 0 & \cdots & 0 \end{pmatrix},
\qquad
\psi_2=\begin{pmatrix} 1 & 0 & -\theta_{13} & \cdots 
& -\theta_{1n} \\ 0 & 1  & -\theta_{23} & \cdots 
& -\theta_{2n} \\ -\theta_{31} & -\theta_{32}  & 0 & \cdots & 0 \\
\cdots & \cdots  & \cdots & \cdots & \cdots \\
-\theta_{n1} & -\theta_{n2}  & 0 & \cdots & 0 \end{pmatrix},$$
and for $3\le k\le n$, 
$$\psi_k=\begin{pmatrix} 0 & 0 & \cdots &\cdots &  \theta_{1k} & \cdots & \cdots \\
0 & 0 & \cdots &\cdots &  \theta_{2k} &\cdots &  \cdots \\
\cdots & \cdots  & \cdots & \cdots & \theta_{3k} & \cdots & \cdots \\
\cdots & \cdots  & \cdots & \cdots & \cdots & \cdots & \cdots \\
\theta_{k1} & \theta_{k2} & \theta_{k3}  & \cdots & 1 & \cdots & 
-\theta_{kn} \\
\cdots & \cdots  & \cdots & \cdots & \cdots & \cdots & \cdots \\
\cdots & \cdots  & \cdots & \cdots & -\theta_{nk} & \cdots & \cdots \\
 \end{pmatrix}.$$
Indeed, the reader can check that for each $t$, the commutators of these
matrices have order two with respect to the local parameters $\mu,\nu$,
and $\theta_{kl}$. 

This defines a basis of $T_{\fa(t)}Red(n)$, by associating to each 
local parameter the tangent vector in the corresponding direction. 
For example, to the parameter $\mu$, we associate the homomorphism
$\partial/\partial\mu \in Hom(\fa,\fsl_n/\fa)$ mapping $\psi_k(t,0)$
to $\partial\psi_k/\partial\mu (t,0)$. Explicitly:
\begin{eqnarray}\nonumber
 &\partial/\partial\mu(\psi_1(t,0)) =e_1^*\ot e_1-e_2^*\ot e_2, & 
\partial/\partial\mu(\psi_k(t,0))=0 \, \,{\rm for}\,k>1, \\
\nonumber
  &\partial/\partial\nu(\psi_1(t,0)) =e_1^*\ot e_2, & 
\partial/\partial\nu(\psi_k(t,0))=0 \, \, {\rm for}\,k>1,
\end{eqnarray}
and so on. Now, what we have to do is to compare this basis with the
other basis defined by the adjoint action of a Killing orthonormal basis
to $\fa(t)^{\perp}$. For $t\ne 0$, let $t=\tau^2$. Then $\fa(t)\in
Red(n)^0$ is the diagonal algebra associated with the basis 
$e_1+\tau e_2$, $e_1-\tau e_2$, $e_3,\ldots e_n$ of $\CC^n$. Its 
Killing orthogonal has a basis given by $e_1^*\ot e_1-e_2^*\ot e_2$,
$\tau e_2^*\ot e_1-\tau^{-1} e_1^*\ot e_2$ and the $e_j^*\ot e_k$, 
with $j\ne k$ and $j$ or $k$ is bigger than two. This basis is not quite
orthonormal, but the norm of the corresponding volume form 
does not depend on $t$. 

We claim that $\partial/\partial\theta_{jk}=ad(e_j^*\ot e_k)$, as
the reader can check. Moreover,
$$\partial/\partial\mu =\frac{1}{2 \tau}
ad(\tau e_2^*\ot e_1-\tau^{-1} e_1^*\ot e_2), \qquad
\partial/\partial\nu =\frac{1}{4t}ad(e_1^*\ot e_1-e_2^*\ot e_2).$$
Note the factor $\tau$, in agreement with the fact that only the square
of the canonical sheaf is trivial on the open orbit. We deduce that
the squared volume form $\om$ at $\fa(t)$ behaves like
$$\om_{\fa(t)}\simeq (Ct\tau)^2\om_0=C^2t^3\om_0
\qquad {\rm when}\quad t\ra 0,$$
if $\om_0$ denotes the local section of the square of the canonical
bundle defined by our local trivialization. 
Hence a zero of order three along $\cO_{bound}$. Since
the codimension one orbit is itself a quadric section of $Red(n)$, we 
deduce:

\begin{theo} 
The canonical sheaf of the smooth locus $Red(n)_{reg}$ is 
$K_{Red(n)_{reg}}=\cO_{Red(n)_{reg}}(-3)$, up to two-torsion.
\end{theo}

To assert that the canonical sheaf of $Red(n)$ is really 
$\cO_{Red(n)}(-3)$, we would first need to answer the following basic
questions.

\medskip\noindent {\it Question E}. Is $Red(n)$ normal ?  

\medskip\noindent {\it Question F}. Is the Picard group of $Red(n)$
torsion free ? What is its rank ? Is it generated by the hyperplane 
divisor, at least up to torsion ? 

\medskip Note that the hyperplane divisor on $Red(n)_{reg}$ is not
divisible, since $Red(n)_{reg}$ contains lines and even planes,
see Proposition 12.

\subsection{Singularities}

We devote this section to a local study of $Ab(n)$ and $Red(n)$ around 
the closed orbit $\cO''_{min}$. We choose the point of $\cO''_{min}$ 
defined as the space of matrices  whose kernel contains and whose image 
is contained in 
the hyperplane $U=\langle e_1,\ldots ,e_{n-1}\rangle$.
Locally around that point, an $(n-1)$-dimensional subspace of $\fsl_n$ 
is made of matrices of the form
$$\begin{pmatrix} A(u) & e_n^*\ot u \\ \a(u)\ot e_n & -a(u) \end{pmatrix},$$
where $u$ belongs to the hyperplane $U$, $\a$ is a linear 
map from $U$ to $U^*$, $A$ a linear map from $U$ to $End(U)$, and $a =
\trace A$.  
This defines an abelian subalgebra of $\fsl_n$ if and only if the following identities hold:
\begin{eqnarray}
\langle \a(u),v\rangle &= &\langle \a(v),u\rangle , \\
A(v)u-a(u)v &= &A(u)v-a(v)u, \\
\, [A(u),A(v)]w &= &\langle\a(u),w\rangle\ot v-\langle\a(v),w\rangle\ot u, \\
 \langle\a(v),A(u)w\rangle-a(v)\langle\a(u),w\rangle &= 
 &\langle\a(u),A(v)w\rangle-a(u)\langle\a(v),w\rangle .
\end{eqnarray}
Letting $B=A+aI$, we can rewrite these identities as
\begin{eqnarray}
\langle \a(u),v\rangle &= &\langle \a(v),u\rangle , \\
B(u)v &= &B(v)u, \\
\, [B(u),B(v)] &= &\a(u)\ot v-\a(v)\ot u, \\
B(u)^t\a(v) &= &B(v)^t\a(u),
\end{eqnarray}
where $B(u)^t$ is the transpose of $B(u)$ and acts on $U^*$. 

At first order, the third set of equations reduces to $\a(u)\ot v=\a(v)\ot u$ 
and implies that $\a=0$. The second set of equations means that $B$ is mapped to zero
by the map $Hom(U, End(U))=U^*\ot U^*\ot U\ra\Lambda^2U^*\ot U$. We thus get 
$(n-1)^2+(n-1)^2(n-2)/2=n(n-1)^2/2$ independent 
linear equations for the Zariski tangent 
space. Since this is half the dimension of the ambient Grassmannian, the Zariski 
tangent space has dimension $n(n-1)^2/2$, which is bigger than $n(n-1)$ as soon as $n>3$.  
Therefore: 

\begin{prop} 
The minimal orbits $\cO_{min}'$ and  $\cO''_{min}$ 
are contained in the singular locus of $Ab(n)$ for $n\ge 4$. 
\end{prop}

Denote by $A_n$ the projectivized tangent cone to a 
normal slice of $\cO_{min}''$
in $Ab(n)$. 
Equations of this tangent cone are $\a =0$ and the symmetry condition (6) 
on $B$ (these equations define the tangent space), plus the quadratic 
equations implied by (7):
\begin{eqnarray}
u\we v\we  [B(u),B(v)](w)= 0 \qquad \forall u,v,w\in U. 
\end{eqnarray}

Note that the tangent space is parametrized by the 
space of morphisms $B\in Hom(U,End(U))=U^*\ot U^*\ot U$ satisfying 
the symmetry condition (6), which implies that in fact they belong to 
the subspace $S^2U^*\ot U$. This tensor product is the direct sum 
of two irreducible components, $S_{10...0-2}U$ and $U^*$. This copy of 
$U^*$ in the tangent space must correspond to the tangent directions 
to the singular strata isomorphic to $\cO_{min}''\simeq\PP^{n-1}$. 
Since this strata is homogeneous, it is natural to restrict to the 
normal slice given by 
$S_{10...0-2}U$, and characterized by the property that the trace of 
$B$ is identically zero. 
Therefore, $A_n$ is the subvariety of $\PP S_{10...0-2}U$,
defined by the quadratic equations (9). 

\medskip
Let $J$ denote the endomorphism of $\CC^n$ defined by $Je_i=e_{i+1}$, 
where the indexes of the basis vectors are taken modulo $n$. Let $\iota$ 
denote the inclusion of $U$ in $\CC^n$, and  $\pi$ the 
projection to $U$ along $e_n$. Let $B(e_i)=\pi J^i\iota$. 
We first claim that $B$
belongs to $A_n$. Indeed, we have 
$$B(e_i)B(e_j)(v)=B(e_i)\sum_{k+j\neq n}v_ke_{j+k}=\sum_{j+k,i+j+k\neq n}v_ke_{i+j+k},$$
so that the commutator $[B(e_i),B(e_j)]$ is simply given by 
$$[B(e_i),B(e_j)](w)=w_{n-i}e_j-w_{n-j}e_i.$$
We immediately deduce that 
$$[B(u),B(v)](w)=(\sum_iw_{n-i}u_i)v-(\sum_jw_{n-j}v_j)u,$$
and the equations (9) follow. 
 
We can be a little more precise: $B$ defines a tangent direction not
only to $Ab(n)$, but really to $Red(n)$. This is because
the space of matrices 
$$\fa(t)=\Bigg\{
\begin{pmatrix} 0 & tx_{n-1} & \cdots & tx_2 & x_1 \\
tx_1 & 0 & \cdots & tx_3 & x_2 \\
\cdots &  \cdots &  \cdots &  \cdots & \cdots \\
tx_{n-2} & \cdots & tx_1  & 0 & x_{n-1} \\  
t^2x_{n-1} & \cdots & t^2x_2 & t^2x_1  & 0 
\end{pmatrix} \Bigg\}$$
is an abelian $(n-1)$-dimensional subalgebra of $\fsl_n$, 
passing through our prefered point of $\cO_{min}''$, whose 
generic point is in $Red(n)^0$ (since, for example, if we let 
$x_i=0$ for $i>0$ we obtain a regular semisimple matrix when 
$t^n\neq 1$). We thus 
get a rational curve on $Red(n)$ whose tangent direction at
$t=0$ is precisely defined by $B$. 

\begin{lemm}
The stabilizer $K_n$ of $B$ in $PGL_{n-1}$ is finite. 
\end{lemm}

\proof The Lie algebra of this stabilizer 
is the space of endomorphisms $X\in\fsl_n$ such that 
$$X.B(u)(v)+B(Xu)(v)+B(Xv)(u)=0 \qquad \forall u,v\in U.$$
Let $f_i=Xe_i$, and $f_i^{+j}=B(e_j)f_i$. We get the conditions 
$$f_{i+j}+f_i^{+j}+f_j^{+i}=0,$$
where indices are taken modulo $n$ and with the convention that $f_n=0$.
We deduce  that $f_2=-2f_1^{+1}$, then $f_3=2f_1^{+1+1}-f_1^{+2}$. 
Then the condition that $f_1^{+3}+f_3^{+1}=2f_2^{+2}$ implies that 
$f_1$ is a combination of $e_{n-1}$
and $e_{n-3}$, and the condition that
$f_1^{+4}+f_4^{+1}=f_2^{+3}+f_3^{+2}$ leads to 
$f_1=0$. Then by induction $f_i=0$ for all $i$, that is, $X=0$. \qed

\medskip\noindent {\it Question G}. What is this finite group $K_n$? 

\medskip A consequence of the lemma is 
that the orbit of the tangent direction defined by $B$ 
has dimension $(n-1)^2-1$, 
which is exactly one minus the dimension $n(n-1)-(n-1)$ of 
our normal slice to $\cO_{min}''$ in $Red(n)$. 
This suggests the following question:
 
\medskip\noindent {\it Question H}. Does the closure $C_n$ of this 
orbit coincide with the projectivized tangent cone to the normal
slice to $\cO_{min}''$ in $Red(n)$? 

\medskip\noindent {\it Question I}. 
When is $C_n$ a smooth variety ? What can its singularities be? 

\medskip We cannot say much about this compactification 
$C_n$ of $PGL_{n-1}/K_n$, which should be an interesting 
object of study. We'll prove in the next section that $C_4$
is in fact a familiar object. 
At least can we say that for $n\ge 4$, $C_n$ is not a linear space
-- and we can therefore conclude:

\begin{prop} 
The minimal orbits $\cO_{min}'$ and  $\cO''_{min}$ 
are contained in the singular locus of $Red(n)$ for $n\ge 4$. 
\end{prop}

\noindent {\it Question J}. When does $Red(n)_{sing}= 
\cO'_{min}\coprod \cO''_{min}$ ?

\subsection{Linear spaces and the incidence variety}

For $n=3$ we proved in \cite{imcrelle} that through a general point of 
$Red(3)$, 
there passes exactly three 
projective planes, which are transverse, and maximal. 
(In fact $Red(3)$ does not contain 
any linear space of dimension greater than two.) 
This extends to $Red(n)$ for any $n$:

\begin{prop}
Through a general point of $Red(n)$, there passes $\binom{n}{2}$ 
projective planes, which are
transverse, and maximal linear subspaces of $Red(n)$. 
\end{prop}

\proof A general point of $Red(n)$ is an  $n$-plane $E$ of $\fgl_n$ 
generated by 
$e_1^*\ot e_1, \ldots, e_n^*\ot e_n$ for some basis $e_1,\ldots ,e_n$ 
of $\CC^n$ and its
dual basis $e_1^*,\ldots ,e_n^*$. A linear space in $G(n-1,\fsl_n)$
 passing through $pE$
is defined by two subspaces $P\subset pE\subset Q$ of $\fsl_n$, where 
$P$ is a hyperplane 
in $pE$ or $pE$ a hyperplane in $Q$. 

In the first case, we get a linear space contained in $Red(n)$ if and 
only if $Q$ is contained
in the centralizer of $P$. Thus $P$ cannot contain any regular element, 
since otherwise its 
centralizer would be equal to $E$. Therefore $P$ must be defined by the 
condition that two 
vectors $e_i$ and $e_j$ belong to the same eigenspace, and its
centralizer 
in $\fsl_n$ has 
dimension $n+1$. We thus get $\binom{n}{2}$ linear spaces in $Red(n)$ 
through $pE$, which are 
all projective planes, and clearly transverse. 

In the second case, $Q$ is generated by $pE$ and some non diagonal 
endomorphism $x$,
which we can suppose to have zero diagonal coeficients. 
A hyperplane in $Q$ not containing $x$ is the space of endomorphisms 
of the form 
$t+\mu (t)x$, with $t\in E$, for some linear form $\mu$ on $E$. It will 
be an abelian
subalgebra of $\fsl_n$ if and only if $\mu(t)[x,s]=\mu(s)[x,t]$ for all 
$s,t\in E$. 
This can hold only if $[x,t]=\mu (t)y$ for some $y\in\fsl_n$,  which 
implies that 
$x=e_i^*\ot e_j$ for some $i\neq j$. But then $\mu$ is fixed up to 
constant, so our linear
space is at most a line and we are back to the first case. \qed   

\medskip
Using this fact we can investigate  the structure of $Red(n)$ through 
the incidence variety
$Z_n$ defined by the diagram
$$\begin{array}{ccccc}
 & & Z_n & & \\
 & \swarrow & & \searrow & \\
\PP\fsl_n & &\dashrightarrow  & & Red(n)\subset G(n-1,\fsl_n)
\end{array}$$
The map $\pi :Z_n\ra Red(n)$ is a $\PP^{n-2}$-bundle, the restriction 
to $Red(n)$ of the tautological vector bundle over $G(n-1,\fsl_n)$. 
The projection $\sigma : Z_n\ra\PP\fsl_n$ is birational. 

\begin{prop}
The map $\sigma : Z_n\ra\PP\fsl_n$ is an isomorphism exactly above the 
open subset 
 of regular elements of $\fsl_n$. 
\end{prop}

Recall that an endomorphism is {\em regular} if its centralizer has
minimal 
dimension. This 
means that there is only one Jordan block for each eigenvalue. The set 
$\bar{W}_n$ of non
regular elements is irreducible, with an open subset given by the set of 
semi-simple 
elements with a double eigenvalue. In particular, $\bar{W}_n$ is the 
projection to $\fsl_n$,
of the set  $W_n\subset\fgl_n$ of elements of corank at least two.   

\medskip
\proof The fiber of $\s$ over $x\in\fsl_n$ is the space of 
$(n-1)$-dimensional abelian 
subalgebras $\fa\subset\fsl_n$ containing $x$, hence contained the 
centralizer 
$\fc_0(x)=\fc(x)\cap\fsl_n$. If $x$ is regular, $\fc(x)$ has dimension 
$n$, hence 
$\fa=\fc_0(x)$, so that $\s$ is one-to-one
over $x$. 

Now suppose that $x$ is a generic non regular endomorphism, so that $x$ 
is semisimple
with distinct eigenvalues, except one with multiplicity two. Let 
$P\simeq\PP^1$ 
denote the projective line defined by its two-dimensional eigenspace. 
Our abelian subalgebra 
$\fa$ is defined by an abelian two-dimensional subalgebra of $\fgl(P)$, 
which is either 
a maximal torus defined by a pair of distinct points in $P$, or the 
centralizer of a nilpotent 
element, defined by one point in $P$. We conclude that the fiber of 
$\s$ over $x$ 
is $Sym^2\PP^1\simeq\PP^2$. Therefore $\s$ in not one to one over the 
generic point 
of $x$, and a fortiori over the whole of $\bar{W}_n$.\qed

\medskip Note that $\bar{W}_n$ has codimension three. Since the fiber 
of $\s$ over its generic 
point has dimension two, we get an exceptional divisor $E\subset Z_n$ 
dominating $\bar{W}_n$,
whose generic point is a pair $(x\in\ft)$, with $\ft$ a Cartan
subalgebra of $\fsl_n$ and $x$ a non 
regular element of $\ft$. The intersection of $E$ with the generic fiber 
$p^{-1}(\ft)$ of $p$,
where $\ft$ denotes the diagonal torus, is the union of the
$\binom{n}{2}$ 
hyperplanes 
$H_{ij}\subset\PP\fsl_n$ defined by the equations $t_i=t_j$, where $i<j$. 

Note that the intersections of these hyperplanes are of two different
types: points in the intersections $H_{ij}\cap H_{kl}$ map to points
in $\fsl_n$ such that the line from that point to the identity matrix
$I$ is bisecant to $W_n$;  points in the intersections $H_{ij}\cap
H_{jk}$ map to the projection to $\PP\fsl_n$ of the variety of 
 matrices of corank at least three.

Since $\sigma$ is birational, we get an induced rational map
$\varphi :\PP\fsl_n\dashrightarrow  Red(n)$ of relative dimension 
$n-2$ which has quite interesting properties. 
First note that 
if $\ell$ is the class of a line in a general fiber of $\pi$, we have  
$E.\ell=\binom{n}{2}$. But $H.\ell=1$, hence $\pi^*\cO(1)=\binom{n}{2}H-E$. 
This proves:

\begin{prop}
The rational map $\varphi :\PP\fsl_n\ra Red(n)$ is defined by a 
linear system $I_n$ of 
polynomials of degree $\binom{n}{2}$, whose base locus contains  $\bar{W}_n$.
\end{prop}

Therefore, we need to understand  $\bar{W}_n$ a little better. Geometrically,
we have the following simple description. 

\begin{prop}
The variety  $\bar{W}_n$ is the projection from the identity matrix,
of the variety of matrices  of corank at least two in $\PP\fgl_n$. 
\end{prop}

\proof By definition, a matrix $X\in\fsl_n$ is non regular if and only if 
some eigenvalue $\lambda$ has multiplicity $m\ge 2$. But then $X-\lambda
I$ has corank $m$ and projects to $x$. The converse assertion is not less
obvious. \qed

\medskip Note that if $X\in\fsl_n$ has an
eigenvalue with  multiplicity two or more, 
its minimal polynomial has degree less than $n$ -- and 
conversely. We deduce:

\begin{prop}
A matrix $X\in\fsl_n$ belongs to the cone over  $\bar{W}_n$ 
if and only if $X,pX^2,\ldots,pX^{n-1}$ are linearly dependent.
\end{prop}

\noindent {\it Question K}. Does this condition define $\bar{W}_n$ 
scheme-theoretically ? 

\medskip
For sure there is no equation of $\bar{W}_n$ of degree smaller than
$\binom{n}{2}$, since the intersection of $\bar{W}_n$ with any Cartan 
subalgebra of $\fsl_n$ is a collection of $\binom{n}{2}$ hyperplanes. 

\smallskip The previous proposition motivates the introduction of a
map 
$$\begin{array}{cccc}
t_n : & \Lambda^{n-1}\fsl_n & \ra & S^{\binom{n}{2}}\fsl_n \\
      &    \theta  & \mapsto & \theta (X,pX^2,\ldots ,pX^{n-1}),
\end{array}$$
where we recall that 
$p : \fgl_n\ra\fsl_n$ denotes the natural projection, and  we 
identify $\fsl_n$ with its dual through the trace map. 
More explicitely, if 
$Z_1,\ldots,Z_{n-1}\in\fsl_n$, we have 
$$t_n(Z_1\cwedge Z_{n-1})(X)=\det (\trace (Z_iX^j))_{1\le i,j\le n-1}.$$
As we have just seen, the base locus of
the image of $t_n$ is equal to $\bar{W}_n$.

The exterior power $\Lambda^{n-1}\fsl_n$ can in principle be decomposed 
as a $GL_n$-module as follows: in the Grothendieck ring of finite
dimensional $GL_n$-modules, we have the identity
$$\Lambda^{n-1}\fsl_n=\sum_{k=1}^n(-1)^{k-1}\Lambda^{n-k}\fgl_n.$$
These wedge powers, since $\fgl_n=\CC^n\ot (\CC^n)^*$, can be decomposed 
using the Cauchy formulas, and then the Littlewood-Richardson rule can
be used to perform the tensor products. This is not very 
effective, and in fact the problem
of decomposing the exterior powers of the adjoint representation 
of $\fsl_n$, and more generally of a simple Lie algebra, has been
much studied since the pionnering work of B. Kostant \cite{kos}, 
see also \cite{bz,reeder} and references therein. 

\smallskip
Of course the map $t_n$ above has no reason of being injective -- we'll see 
in the next section that injectivity fails already for $n=4$. 
This leaves quite a number of open questions:

\medskip\noindent {\it Question L}. Does
$I_n=I_{\bar{W}_n}(\binom{n}{2})$ ?
Does $I_n=Ker \;\Theta$ ? 
Does $I_n={\rm Im}\; t_n$ ? And can
we compute ${\rm Im}\; t_n$ explicitly ? 

\medskip
A simple fact to mention about the image
of $t_n$ is that it certainly contains $S^n\CC^n$ and its dual, embedded 
through the map
$$\begin{array}{cccc}
s_n : & S^n\CC^n\ot (\det \CC^n)^{-1} & \ra & S^{\binom{n}{2}}\fsl_n \\
      &    v^n\ot (u_1\cwedge u_n)^{-1} & \mapsto & P_v(X)=
v\we Xv\cwedge X^{n-1}v/u_1\cwedge u_n,
\end{array}$$
and similarly, the dual map 
$$\begin{array}{cccc}
s_n^* : & S^n(\CC^n)^*\ot (\det \CC^n) & \ra & S^{\binom{n}{2}}\fsl_n \\
      &    e^n\ot (f_1\cwedge f_n)^{-1} & 
\mapsto & P^*_e(X)=e\we { }^t Xe\cwedge { }^t X^{n-1}e/f_1\cwedge f_n.
\end{array}$$
(It is enough to define these maps on pure powers, since they generate
the space of polynomials. The image of a monomial can be deduced by
polarization.)

\smallskip
For future use, note that we have a commutative diagram
$$\begin{array}{ccccc}$$
S^n\CC^n & \stackrel{i}{\lra} & \Lambda^{n-1}\fsl_n &
\stackrel{\pi}{\lra} & S^n\CC^n \\
 || & & \stackrel{\a}{}\downarrow\uparrow \stackrel{\b}{} & & || \\
S^n\CC^n & \stackrel{j}{\lra} & \Lambda^n\fgl_n &
\stackrel{\rho}{\lra} & S^n\CC^n,
\end{array}$$
defined as follows. First, the map $j$ is given by 
$$ j(v^n\ot e_1^*\cwedge e_n^*)=(e_1^*\ot v)\cwedge (e_n^*\ot v).$$
Second, we define the projection $\rho$ by letting
$$\rho (X_1\cwedge X_n)=\det (X_i(e_j))\ot e_1^*\cwedge e_n^*,$$
where $e_1,\ldots ,e_n$ is any basis of $\CC^n$ and $e_1^*,\ldots ,e_n^*$
the dual basis. Note that $\rho\circ j=id$.

Of course the maps $\a$ and $\b$ are defined through the decomposition 
$\Lambda^n\fgl_n=\Lambda^n(\CC I\op\fsl_n)=\Lambda^{n-1}\fsl_n\op
\Lambda^n\fsl_n$. Explicitly, we have 
\begin{eqnarray}
\nonumber
\a (Y_1\cwedge Y_{n-1}) & = & I\wedge Y_1\cwedge Y_{n-1}, \\
\nonumber
\b (X_1\cwedge X_n) & = & \frac{1}{n}\sum_{j=1}^n(-1)^{j-1}\tr (X_j)X_1\cwedge
\hat{X_j}\cwedge X_n.
\end{eqnarray}
The choice of the constant $\frac{1}{n}$ 
is such that $\b\circ\a=id$. Finally, we let 
$i=\b\circ j$ and $\pi=\rho\circ\a$. We have $\pi\circ i=\frac{1}{n}id$,
since
\begin{eqnarray} 
\nonumber
v^n\ot e_1^*\cwedge e_n^* & \stackrel{i}{\mapsto} & \frac{1}{n}
\sum_j (-1)^{j-1}\langle e_j^*,v\rangle 
(e_1^*\ot v)\cwedge (e_{j-1}^*\ot v)\cwedge (e_{j+1}^*\ot v) \\
\nonumber
 & \stackrel{\a}{\mapsto} &  \frac{1}{n}
\sum_j \langle e_j^*,v\rangle 
(e_1^*\ot v)\cwedge I\cwedge (e_n^*\ot v) \\
\nonumber
 & \stackrel{\rho}{\mapsto} & \frac{1}{n}
\sum_j \langle e_j^*,v\rangle \det\begin{pmatrix} v & & e_1&  & \\
 & v & & & \\ & & e_j & & \\ & & & v & \\ & & e_n & & v\end{pmatrix}
\ot e_1^*\cwedge e_n^* \\
\nonumber &  = & \frac{1}{n}\sum_j \langle e_j^*,v\rangle  v^{n-1}e_j\ot 
e_1^*\cwedge e_n^* =\frac{1}{n}v^n\ot 
e_1^*\cwedge e_n^*.
\end{eqnarray} 

We use this diagram to define an automorphism $\tau$ of
$\Lambda^{n-1}\fsl_n$ by 
$$\tau=id-(-1)^{n-1}\frac{(n-1)!}{n}\; i\circ\pi,$$
and we twist our map $t_n$ above by letting $t'_n=t_n\circ\tau$.
The point is that we now have:

\begin{prop}
Let $\ft\in\Lambda^{n-1}\fsl_n$ belong to the cone over $Red(n)$. 
Then the polynomial $t'_n(\ft)$ on $\fsl_n$ vanishes on the linear
subspace $\ft$.
\end{prop}

\proof We just need to prove it for a generic element of $Red(n)$, 
that is, we may suppose that $\ft$ corresponds to the space of matrices
that are diagonal with respect to some basis $e_1,\ldots ,e_n$. 
Up to constant, we may therefore let
$$\ft = (e_1^*\ot e_1-e_2^*\ot e_2)
\cwedge (e_{n-1}^*\ot e_{n-1}-e_n^*\ot e_n).$$
If $x\in\ft$ is diagonal with eigenvalues $x_1,\dots ,x_n$, we 
immediately get 
$$t_n(\ft)(x)=\det (x_i^j-x_{i+1}^j)=(-1)^{n-1}\prod_{i<j}(x_i-x_j).$$
On the other hand, let us compute $\ft_0=i\circ\pi(\ft)$. First note
that
$$\ft=\sum_{j=1}^{n}(-1)^{j-1}(e_1^*\ot e_1)\cwedge (e_{j-1}^*\ot
e_{j-1})\wedge (e_{j+1}^*\ot e_{j+1})\cwedge (e_n^*\ot e_n),$$
and since of course $I=e_1^*\ot e_1+ \cdots +e_n^*\ot e_n$, we deduce
that $\a(\ft)=n(e_{1}^*\ot e_{1})\cwedge (e_n^*\ot e_n)$, hence 
$\pi (\ft)=n e_1\cdots e_n\ot e_1^*\cwedge e_n^*$ and 
$$j\circ\pi (\ft)=\frac{1}{(n-1)!}\sum_{\s\in\cS_n}\e(\s) 
(e_{\s(1)}^*\ot e_{1})\cwedge (e_{\s(n)}^*\ot e_n).$$
From that we could easily compute $\ft_0$, but to finish the computation
we prefer to notice that for all $Y_1,\ldots ,Y_{n-1},X\in\fgl_n$, the 
identity 
$$(pY_1\cwedge pY_{n-1})(pX,\ldots ,pX^{n-1})=(I\wedge Y_1\cwedge Y_n)
(I,X,\ldots ,X^{n-1})$$
holds true, as the reader can easily check. Therefore,
\begin{eqnarray}
\nonumber
\b(Y_1\cwedge Y_n)(pX,\ldots ,pX^{n-1}) & = & \sum_{j=1}^{n}(-1)^{j-1}
(Y_1\cwedge \hat{Y_j}\cwedge Y_n) (pX,\cdots ,pX^{n-1}) \\
\nonumber {  } &= & \sum_{j=1}^{n}(Y_1\cwedge Y_{j-1}\wedge I\wedge Y_{j+1}
\cwedge Y_n) (I,X,\ldots ,X^{n-1}).
\end{eqnarray}
Applying this to $X\in\ft$ a diagonal matrix and $Y_i=e_{\s(i)}^*\ot
e_i$, we see that only $\s=1$ can contribute. Then we get $I=Y_1+
\cdots+Y_n$, and we finally deduce that
$$t_n(\ft_0)(X)=\frac{n}{(n-1)!}\prod_{i<j}(x_i-x_j).$$
This completes the proof. \qed

\medskip If, as we expect, $I_n=Im\; t_n= Im\; t'_n$, this proposition gives a 
coherent identification between $Red(n)$ and the image of the rational
map defined by that linear system. Indeed, the image of a general point 
$x$ is  the hyperplane of equations of ${\bar W}_n$ also
vanishing on $x$. By the proposition, $s_n(x)$, seen as a
hyperplane of equations, does vanish on the centralizer of $x$, in 
particular on $x$ itself.  
Without the twist $\tau$, the image of the rational map defined by $I_n$
is only a translate of $Red(n)$ by a linear automorphism.

\subsection{Relations with the Hilbert scheme}

The open $PGL_n$-orbit of $Red(n)$ is the space of $n$-tuples of points
in general
position in $\PP^{n-1}$. This is also an open $PGL_n$-orbit in the
punctual Hilbert 
scheme $Hilb^n\PP^{n-1}$, which is thus birational to $Red(n)$. For $n=3$ we 
proved in \cite{imcrelle} that 
there is a morphism $Hilb^3\PP^2\ra Red(3)$, in fact a divisorial
contraction 
between
these two smooth varieties. We would like to extend this to $n\ge 4$. 

We define auxiliary morphisms as follows. First, we have a 
$GL_n$-equivariant map 
$$\begin{array}{cccc}
\mu_n : &  \Lambda^n(S^{n-1}\CC^n) & \lra & \Lambda^n(S^{n-1}\CC^n) \\
 & e_1^{n-1} \cwedge e_n^{n-1} & \mapsto & e_2\cdots e_n\cwedge
  e_1\cdots 
e_{n-1}.
\end{array}$$

\medskip\noindent {\it Question M}. Is $\mu_n$ an isomorphism for all
$n$ ? (We know it is for $n=3$.)  

\medskip Now define the $SL_n$-equivariant morphism
$$\begin{array}{cccc}
\nu_n : &  \Lambda^n(S^{n-1}\CC^n) & \lra & \Lambda^n\fgl_n \\
 & e_1^{n-1} \cwedge e_n^{n-1} & \mapsto & (e_1\cwedge e_{n-1}\ot e_n)
\cwedge (e_n\wedge e_2\cwedge e_{n-1}\ot e_1).
\end{array}$$
Here we identify $\Lambda^{n-1}\CC^n$ with $(\CC^n)^*$, hence  
$\Lambda^{n-1}\CC^n\ot\CC^n$ with $\fgl_n$. 
For example, if $e_1,\ldots, e_n$ are independent and $e_1^*,
\ldots ,e_n^*$ is the dual basis, 
the image tensor is just $(e_1^*\ot e_1)\cwedge (e_n^*\ot e_n)$. 
This can be identified (up to scalar, 
of course), with the linear space generated by the $n$ rank one 
elements $e_k^*\ot e_k$ of $\fgl_n$,
the sum of which is the identity. In other words, we get a point of the 
open orbit of our reduction 
variety $Red(n)$. 

\medskip Let $Z$ be a $n$-tuple of points in general position in
$\PP^{n-1}$. 
Denote by $T(Z)$ 
the union of the $\binom{n}{2}$ codimension two linear spaces generated
by 
all the $(n-2)$-tuples 
of points in $Z$. Let $P$ denote the Hilbert polynomial of this
variety $T(Z)$. Once we have chosen homogeneous coordinates 
adapted to our $n$-tuple of points, we see that the ideal of $T(Z)$
is generated by the monomials $x_1\cdots {\hat x_i}\cdots {\hat x_j}
\cdots x_n$, so that a supplement of $I_{T(Z)}(k)$ has a basis 
given by all the degree $k$ monomials involving no more than 
$n-2$ indeterminates. We deduce that 
$$P(k)=\sum_{j=1}^{n-2}\binom{n}{j}\binom{k-1}{j-1}.$$
The transformation $T$ defines a rational map
from $Hilb^n\PP^{n-1}$ to $Hilb^P\PP^{n-1}$, which is 
not a morphism in general.  
But we may be able to define a morphism $\rho_n : 
Hilb^n\PP^{n-1}\ra Red(n)$ 
as follows. We first map the punctual scheme $Z$, reduced and in general
position, to the linear system 
$|I_{T(Z)}(n-1)|\in
G(n,S^{n-1}\CC^n)\subset\PP\Lambda^n(S^{n-1}\CC^n)$. Then we apply the 
linear automorphism 
$\mu_n^{-1}$, and finally the linear morphism $\nu_n$ to get a point 
$\rho_n(Z)\in
\PP\Lambda^n\fgl_n$. We claim that $\rho_n$ maps the component 
$Hilb_0^n\PP^{n-1}$ of 
$Hilb^n\PP^{n-1}$ containing the reduced schemes 
in general position, to the reduction variety 
$Red(n)$. Indeed, if $Z$ is the union of $n$ points in general position, 
and if $e_1,
\ldots ,e_n$ are adapted coordinates,  then 
$$|I_{T(Z)}(n-1)| = \langle e_1\cdots e_{n-1}, \ldots , e_2\cdots 
e_n\rangle, $$
so that $\rho_n(Z)$ is the subspace of $\fgl_n$ generated by $e_1^*\ot 
e_1,\ldots , e_n^*\ot e_n$,
and belongs to $Red(n)$.

\medskip\noindent {\it Question N}. Can $\rho_n$ be extended
to a morphism ? If yes, 
what is the exceptional locus of  this morphism ? Is $\rho_n$ a 
divisorial contraction for $n\ge 4$ ?

\pagebreak
\section{Reductions for $\fgl_4$}

Its seems rather difficult to answer in full generality  the questions 
we raised in the first part of this paper. In this second part we check
that (almost) everything works as expected when $n=4$. Our first result 
is that $Red(4)=Ab(4)$. More precisely:

\begin{prop}
The variety of reductions in $\fsl_4$, coincides with the space of 
three-dimensional
abelian subalgebras of $\fsl_4$.
It is made of $14$ $PGL_4$-orbits, exactly three of which 
are closed: a three-dimensional projective space $\PP^3$ and its dual 
${\check \PP}^3$, 
and a variety of complete flags $\FF_4$. 
\end{prop}

The proof of this result will occupy the next two sections. 

\subsection{Classification of three dimensional abelian subalgebras 
of $\fsl_4$}

First, we have the one-regular abelian subalgebras, whose different
types are given by the possible sizes of the Jordan blocks of a generic 
element. We thus get five {\it regular} orbits, with generic Jordan type
$1111$ (genuine reductions), $211$, $22$, $31$ or $4$ (regular
nilpotents) (the numbers are just the sizes of the Jordan blocks). 
We denote these orbits by $\cO_{12}$, $\cO_{11}$, $\cO_{10}$, 
$\cO'_{10}$, $\cO_9$ respectively.
 
\smallskip
Now suppose that $\fa$ contains no regular element. If it contains an element 
of Jordan type $211$, $\fa$ is contained in its centralizer which is a copy of 
$\fgl_2\times\fgl_1\times\fgl_1$, and the blocks from $\fgl_2$ are generically 
non regular. But in dimension two this means that they are homotheties, and
this leaves only two free parameters, a contradiction. The Jordan type $22$ 
is eliminated for the same reason. If $\fa$ contains an element of Jordan type
$31$, it must be contained in $\fgl_3\times\fgl_1$ and the blocks from $\fgl_3$
must be non regular, hence of the form $xI+X$ with $X^2=0$ and we need an 
abelian plane of such endomorphisms. We know this leaves only two 
possibilities (in fact only one up to transposition),
$$\begin{pmatrix} c & 0 & a & 0 \\ 0 & c & b & 0 \\
0 & 0 & c & 0 \\ 0 & 0 & 0 & d \end{pmatrix}
\quad {\rm or} \quad 
\begin{pmatrix} c & b & a & 0 \\ 0 & c & 0 & 0 \\
0 & 0 & c & 0 \\ 0 & 0 & 0 & d \end{pmatrix}.$$
Hence two orbits  $\cO'_8$ and  $\cO''_8$.
 
We are left with the {\it nilpotent} abelian algebras containing no
regular element. 
If there is an element $x$ with a Jordan block of size 3, say 
$$x=\begin{pmatrix} 0 & 1 & 0 & 0 \\ 0 & 0 & 1 & 0 \\
0 & 0 & 0 & 0 \\ 0 & 0 & 0 & 0 \end{pmatrix}
\quad {\rm then} \quad 
y=\begin{pmatrix} a & b & c & d \\ 0 & a & b & 0 \\
0 & 0 & a & 0 \\ 0 & 0 & e & f \end{pmatrix}$$
if $y$ commutes with $x$. If $y$ is nilpotent, $a=f=0$. Since $[y,y']=
(de'-d'e)e_3^*\ot e_1$, 
we'll get a three dimensional abelian algebra if we impose a linear
relation 
between
$e$ and $d$. Up to a change of basis, there are only three cases, $d=e$, 
$d=0$, $e=0$, 
the last two being exchanged by transposition. Hence three orbits
$\cO_8$, 
$\cO'_7$ and  $\cO''_7$.

If no element of $\fa$  has a Jordan block of size 3, then $x^2=0$ for
every 
$x$ in $\fa$. 
Suppose that some $x$ has rank two. Every endomorphism commuting with
$x$ 
will preserve 
its kernel, hence be of the form
$$y=\begin{pmatrix} A & B \\ 0  & C\end{pmatrix},
\quad {\rm so\; that} \quad 
y^2=\begin{pmatrix} A^2 & AB+BC \\ 0  & C^2\end{pmatrix}.$$
Using the commutativity condition, we see that $A$ (and $C$) must vanish 
or be proportional 
to a fixed nilpotent matrix when $y$ varies in $\fa$. If $A$ and $C$ are 
both not identically 
zero, we get up to a change of basis 
$$y=\begin{pmatrix} 0 & a & b & c \\ 0 & 0 & 0 & d \\
0 & 0 & 0 & e \\ 0 & 0 & 0 & 0 \end{pmatrix}, \qquad ad+be=0.$$
This means that $d=ze$, $b=-za$ for some scalar $z$. But then a simple 
change of basis
implies that we may suppose that $A$ and $C$ are in fact both
identically 
zero !
This means that there is a plane $P$ such that every element of $\fa$
has 
$P$ in its kernel
and its image in $P$. In fact this defines a four-dimensional abelian 
algebra, of which
$\fa$ is a hyperplane defined by some non zero linear form. This form 
is defined by some
order two matrix, and changing basis gives the usual $GL_2\times
GL_2$-action 
by left
and right multiplication, with the rank as only invariant. We thus get 
two orbits $\cO_7$
(rank two) and $\cO_6$ (rank one). 

Finally, suppose that $C$ is identically zero, but not $A$. Then the 
condition $AB=0$ 
means that the  the image of $B$ is contained in the kernel of $A$, so 
that $\fa$ is 
the space of traceless endomorphisms with image in a given line. 
Symmetrically, if
$A$ is identically zero, but not $C$, then $\fa$ is 
the space of traceless endomorphisms whose kernel contains a given hyperplane. 
These two orbits $\cO'_3$ and $\cO''_3$ are exchanged by transposition,
they are the minimal orbits denoted $\cO'_{min}$ and $\cO''_{min}$
in the first part of the paper.
\medskip 
Apart from $\cO_8$, $\cO'_7$, $\cO''_7$ and $\cO_7$, all the orbits can 
be described 
in terms of geometric datas. For example, $\cO_{12}$ is the variety of 
quadruples of 
independent points in $\PP^3$. A point in  $\cO_{11}$ is determined by 
two points and
a line in general position, plus a point on the line, and so on. These 
orbits can
therefore be described as open subsets of products of partial flag varieties. 

A point in $\cO'_7$ or $\cO''_7$ determines a complete flag in $\PP^3$, 
and these
orbits are $\CC^*$-bundles over the complete flag variety
$\FF_4$. $\cO_7$ 
is an 
affine fibration over the Grassmannian $G(2,4)$, and $\cO_8$ an affine 
fibration 
over the partial flag variety $\FF_{1,3}$. 

\medskip 
Here is the list of the $14$ orbits with a representative for each. 
(We omit the condition 
that the trace must vanish.) The subscript is the dimension.

$$\begin{array}{lccclcc}
\cO_{12} & \begin{pmatrix} a & 0 & 0 & 0 \\0 & b & 0 & 0 \\
0 & 0 & c & 0 \\0 & 0 & 0 & d \end{pmatrix} &  & \qquad & 
\cO_{11} & \begin{pmatrix} a & b & 0 & 0 \\0 & a & 0 & 0 \\
0 & 0 & c & 0 \\0 & 0 & 0 & d \end{pmatrix} & \\
 & & & & & \\
\cO'_{10} & \begin{pmatrix} a & b & 0 & 0 \\0 & a & 0 & 0 \\
0 & 0 & c & d \\0 & 0 & 0 & c \end{pmatrix} & & & 
\cO''_{10} & \begin{pmatrix} a & b & c & 0 \\0 & a & b & 0 \\
0 & 0 & a & 0 \\0 & 0 & 0 & d \end{pmatrix} & \\
 & & & & & \\
\cO_{9} & \begin{pmatrix} 0 & a & b & c \\0 & 0 & a & b \\
0 & 0 & 0 & a \\0 & 0 & 0 & 0 \end{pmatrix} & & \qquad & 
\cO_{8} & \begin{pmatrix} 0 & a & b & c \\0 & 0 & 0 & a \\
0 & 0 & 0 & b \\0 & 0 & 0 & 0 \end{pmatrix} & \\
 & & & & & \\
\cO'_{8} & \begin{pmatrix} c & b & a & 0 \\0 & c & 0 & 0 \\
0 & 0 & c & 0 \\0 & 0 & 0 & d \end{pmatrix} & & \qquad &
\cO''_{8} & \begin{pmatrix} c & 0 & a & 0 \\0 & c & b & 0 \\
0 & 0 & c & 0 \\0 & 0 & 0 & d \end{pmatrix} & \\
 & & & & & \\
\cO'_{7} & \begin{pmatrix} 0 & a & b & c \\0 & 0 & 0 & a \\
0 & 0 & 0 & 0 \\0 & 0 & 0 & 0 \end{pmatrix} & & & 
\cO''_7{} & \begin{pmatrix} 0 & 0 & a & c \\0 & 0 & 0 & b \\
0 & 0 & 0 & a \\0 & 0 & 0 & 0 \end{pmatrix} & \\
 & & & & & \\
\cO_{7} & \begin{pmatrix} 0 & 0 & a & b \\0 & 0 & b & c \\
0 & 0 & 0 & 0 \\0 & 0 & 0 & 0 \end{pmatrix} & & & 
\cO_{6} & \begin{pmatrix} 0 & 0 & a & b \\0 & 0 & 0 & c \\
0 & 0 & 0 & 0 \\0 & 0 & 0 & 0 \end{pmatrix} & \FF_4 \\
 & & & & & \\
\cO'_{3} & \begin{pmatrix} 0 & a & b & c \\0 & 0 & 0 & 0 \\
0 & 0 & 0 & 0 \\0 & 0 & 0 & 0 \end{pmatrix} & \PP^3 & & 
\cO''_{3} & \begin{pmatrix} 0 & 0 & 0 & c \\0 & 0 & 0 & b \\
0 & 0 & 0 & a \\0 & 0 & 0 & 0 \end{pmatrix} & {\Hat \PP}^3
\end{array}$$

\subsection{Degeneracies}

We want to study which orbits are contained in the closure of which. We will
denote $\cO\ra\cO'$ if $\cO'$ is included in the boundary of $\cO$. 

First note that if  $\fa\in\cO$ and $\fa'\in\cO'$ are one-regular, that is, 
can be defined as the centralizers
of some regular elements $x$ and $x'$, we just need to let $x$
degenerate to $x'$
in the open set of regular elements to make $\fa$ degenerate to $\fa'$. And 
letting $x$
degenerate to $x'$ is possible as soon as this is compatible with the
size 
of the 
Jordan blocks. We deduce that $\cO\ra\cO'$ as soon as $\dim\cO >
\dim\cO'$.
More generally, we know that the {\it two-regular} orbits 
in $Ab(4)$ are contained in $Red(4)$. An easy case-by-case 
check leads to the following conclusion:

\begin{lemm}
The only orbits in $Ab(4)$ which are not two-regular are 
$\cO_7$, $\cO_6$, $\cO_3'$ and $\cO_3''$.
\end{lemm}

We complete the picture by showing that any orbit, with of course $\cO_{12}$ 
excepted, is in the closure of an orbit of larger dimension.  This
will imply 
that every three-dimensional abelian subalgebra of 
$\fsl_4$ is contained in the 
closure of the variety of non singular reductions. Actually we prove a little 
more than needed, in order to deduce the full diagram of degeneracies. 

\medskip\noindent 
{\bf $\cO_9\ra\cO_8$}: if we take the representative above of $\cO_9$
and make the change of basis $e_1\ra te_1$, $e_3\ra te_3$, we get the
abelian 
algebra
of matrices of the form
$$\begin{pmatrix} 0 & t^{-1}a & b & t^{-1}c \\0 & 0 & ta & b \\
0 & 0 & 0 & t^{-1}a \\0 & 0 & 0 & 0 \end{pmatrix}=
\begin{pmatrix} 0 & a' & b' & c' \\0 & 0 & t^2a' & b' \\
0 & 0 & 0 & a' \\0 & 0 & 0 & 0 \end{pmatrix},$$
if $a'=t^{-1}a, b'=b, c'=t^{-1}c$. Letting $t\ra 0$, we get an abelian 
subalgebra belonging to $\cO_8$.  

\medskip\noindent 
{\bf $\cO''_{10}\ra\cO'_8$}: if we take the representative above of $\cO''_{10}$
and make the change of basis $e_2\ra t^{-1}e_2$, we get the abelian algebra
of matrices of the form
$$\begin{pmatrix} a & t^{-1}b & c & 0 \\0 & a & tb & 0 \\
0 & 0 & a & 0 \\0 & 0 & 0 & d \end{pmatrix}=
\begin{pmatrix} a & b' & c & 0  \\0 & a & t^2b' & 0 \\
0 & 0 & a & 0 \\0 & 0 & 0 & d \end{pmatrix},$$
if $b'=t^{-1}b$. Letting $t\ra 0$, we get an abelian 
subalgebra belonging to $\cO'_8$. By transposition we also 
have  $\cO''_{10}\ra\cO''_8$.

\medskip\noindent 
{\bf $\cO'_8\ra\cO_7$}: we take the representative above of $\cO'_8$
and make the change of basis $e_1\ra e_4$, $e_2\ra e_3$,
$e_3\ra e_2$,  $e_4\ra e_1+t^{-1}e_4$, and we let $t\ra 0$. 

\medskip\noindent 
{\bf $\cO_8\ra\cO'_7$}: make the change $e_3\ra t^{-1}e_3$ and let $t\ra 0$. 
By transposition we also have  $\cO_8\ra\cO''_7$.

\medskip\noindent 
{\bf $\cO_7\ra\cO_6$}: make the change $e_3\ra t^{-1}e_3$ and let $t\ra 0$
after renormalizing by $c'=tc$. It is not more difficult to see that $\cO'_7\ra\cO_6$
and $\cO''_7\ra\cO_6$. 

\medskip\noindent 
{\bf $\cO'_7\ra\cO'_3$}: make the change $e_2\ra t^{-1}e_2$ and let $t\ra 0$
after renormalizing by $a'=t^{-1}a$. Transposing, we also get 
$\cO''_7\ra\cO''_3$. 

\medskip Finally, we don't have {\bf $\cO_9\ra\cO'_8$} since $\cO_9$ is 
nilpotent but not $\cO_8'$; neither {\bf $\cO_8\ra\cO_7$} because an 
abelian algebra in $\cO_8$ maps a fixed hyperplane to a fixed line, 
while this does not happen for a  abelian algebra in $\cO_7$;
neither {\bf $\cO_7\ra\cO'_3$} or {\bf $\cO_7''\ra\cO'_3$} since the
matrices in an algebra belonging to $\cO'_3$ vanish on a common line but
not on a common plane. 

We deduce the complete  incidence diagram:

$$\begin{array}{cccccc}
 & & \cO_{12}  & & & \\ 
 & &\downarrow  & & & \\ 
 & &\cO_{11}  & & & \\ 
 & \swarrow & & \searrow & & \\ 
\cO'_{10} & & & & \cO''_{10}  & \\
\downarrow & \searrow & & \swarrow &\downarrow & \\  
\downarrow & & \cO_9& & \downarrow  & \\ 
\downarrow &  &\downarrow & & \downarrow & \\  
\cO'_8 & &\cO_8 & &\cO''_8 & \\ 
 & \swarrow\hspace{-4mm}\searrow & & \searrow\hspace{-4mm}\swarrow
 \\
\cO'_7 & &\cO_7 & &\cO''_7 & \\ 
\downarrow & \searrow &\downarrow & \swarrow & \downarrow & \\  
\downarrow & & \cO_6 & & \downarrow  & \\ 
\cO'_3 & & & & \cO''_3 &  
\end{array}$$

\subsection{The linear span of $Red(4)$}

Remember that set-theoretically, $Red(4)=Ab(4)$ can be defined
as a linear section of $G(3,\fsl_4)$ by the kernel of the map
$$\Theta : \Lambda^3\fsl_4\ra\Lambda^2\fsl_4\ot\fsl_4\ra\fsl_4\otimes\fsl_4,$$
obtained by composing the obvious inclusion with the commutator 
$\Lambda^2\fsl_4\ra\fsl_4$. With the help of LiE \cite{LiE}, 
we check that this kernel is
$$ker\Theta = S_{3-1-1-1}\CC^4\op S_{111-3}\CC^4\op S_{21-1-2}\CC^4.$$
Since $Red(4)$ contains three closed orbits $\PP^3$, ${\check\PP}^3$ and 
$\FF_4$
which are the closed orbits in the projectivisations of the simple 
factors of $ker\Theta$, we conclude that
the linear span in $\PP\Lambda^3\fsl_4$ of the abelian subalgebras, 
is the whole of $ker\Theta$. Its dimension is $35+35+175=245$. 

\subsection{The incidence variety and the induced rational map}

Remember the diagram
$$\begin{array}{ccccc}
 & & Z_4 & & \\
 & \swarrow & & \searrow & \\
\PP\fsl_4 & &\dashrightarrow  & & Red(4)\subset G(3,\fsl_4)
\end{array}$$
The map $\pi :Z_4\ra Red(4)$ is a $\PP^{2}$-bundle, while the projection 
$\sigma : Z_4\ra\PP\fsl_4$ is birational, and an isomorphism above the 
open set of regular elements of $\fsl_4$. The rational map $\varphi
 :\PP\fsl_4\dashrightarrow Red(4)$ is defined by a linear system $I_4$
of sextics vanishing on ${\bar W}_4$.

\begin{prop}
The linear system $I_4$ is equal to $I_{{\bar W}_4}(6)$, and to 
the image of $t_4$, and to the kernel of $\Theta$.
\end{prop}

\proof A computation by Macaulay \cite{mac} 
shows that the ideal of ${\bar W}_4$
is generated by $245$ sextics (we thank Marcel Morales for his help
in performing this computation). We already know $245$ such sextics:
the image of $s_4$, a copy of $S^4\CC^4\ot
(\det\CC^4)^{-1}=S_{3-1-1-1}\CC^4$, 
gives $35$ of them; the 
image of $s_4'$ gives $35$ others, a copy of the dual module; 
and the image of $t_4$ contains
$175$ more. Indeed, remember that $t_4$ associates to a triple 
of matrices $Y_1, Y_2, Y_3\in\fsl_4$ the sextic polynomial
$$P(X)=\det \begin{pmatrix} \trace(Y_1X) &\trace(Y_2X) &\trace(Y_3X) \\
\trace(Y_1X^2) &\trace(Y_2X^2) &\trace(Y_3X^2) \\\trace(Y_1X^3) &\trace(Y_2X^3) &\trace(Y_3X^3)
\end{pmatrix}.$$
Choose two independent vectors $u,v\in\CC^4$ and two independent linear forms 
$\a$, $\b$ vanishing on them. Letting $Y_1=\a\ot u$, $Y_2=\b\ot u$ and 
$Y_3=\a\ot v$, we get the polynomial
$$P(X)=\det \begin{pmatrix} \a(Xu) & \b(Xu) & \a(Xv) \\
\a(X^2u) & \b(X^2u) & \a(X^2v) \\\a(X^3u) & \b(X^3u) & \a(X^3v)
\end{pmatrix},$$  
Note that this polynomial remains 
unchanged if we add to $v$ a multiplle of $u$, or to $\b$ a multiple 
of $\a$. This means that, up to constant, this 
polynomial only depends on the complete flag $\CC u\subset\CC u\op\CC v=
Ker(\a)\cap Ker(\b)\subset\ Ker(\a)$. We conclude that the projectivized
image of $t_4$ contains a copy of the compete flag manifold $\FF_4$. 
Moreover, since the weights of $u,v,\b,\a$ in $P$ are $2,1,1,2$, 
the linear span of  this flag manifold is a copy of the $GL_4$-module 
$S_{21-1-2}\CC^4$, which has dimension $175$. We conclude that, as 
a $GL_4$-module, 
$$I_{{\bar W}_4}(6)=S_{3-1-1-1}\CC^4\op S_{111-3}\CC^4\op
S_{21-1-2}\CC^4=Im t_4.$$
This is isomorphic with $Ker\Theta$; more precisely, $t_4$ restricts to 
an isomorphism between $Ker\Theta$ and $Im t_4$, 
since a computation by LiE shows that $\Lambda^3\fsl_4$
is multiplicity free.\qed 

\medskip 
Once this is established, we can understand the map $\varphi$
geometrically, in particular we can describe 
the fiber of $\sigma$ over most points
$x\in\fsl_4$, that is, the variety parametrizing the abelian three-dimensional
subalgebras of $\fsl_4$ cointaining $x$. 

To state our next result we need to define several natural 
subvarieties of $\PP\fsl_4$. We already introduced the variety
${\bar W}_4$ of non regular elements, and the projection ${\bar X}_4$
of the rank one variety in $\PP\fgl_4$. Let $X_4^0$ denote the variety 
of rank one matrices in $\PP\fsl_4$. 

Let also $Y_4$ denote the space of matrices in $\PP\fsl_4$ which belong 
to some bisecant line to $W_4$ passing through the identity matrix. A generic 
point in $Y_4$ is a matrix with two (opposite) eigenvalues, both of
multiplicity two, so that $Y_4$ contains an open $PGL_4$-orbit
isomorphic with the space of pairs of skew lines in $\PP^3$. Let 
$Y^4_0$ denote the complement of the open orbit in $Y_4$. The points
of $Y^4_0$ are nilpotent matrices, either with two Jordan
blocks of size two, or of rank one (hence in $X_4^0$).  The singular 
locus of ${\bar W}_4$ is the union of ${\bar X}_4$ and $Y_4$, whose 
intersection is $X_4^0$. 

\begin{prop}
Let $x\in\PP\fsl_4$. The fiber of $\sigma$ over $x$ is
\begin{enumerate}
\item a point if $x$ is regular, 
\item a projective plane if $x$ belongs to $W_{4,reg}$, 
\item the product of two projective planes if $x\in Y_4-Y_4^0$,
\item a copy of $Red(3)$ if $x\in {\bar X}_4-X_4^0$.
\end{enumerate}
\end{prop}

\proof If $x$ is regular, the unique three-dimensional subalgebra
of $\fsl_4$ that contains it is its centralizer, thus $\s^{-1}(x)$ 
is a point. Note that if moreover $x$ is semisimple, the intersection 
of ${\bar W}_4$ with the Cartan subalgebra $\fc (x)$ is the union 
of six hyperplanes -- so that the linear system $I_{{\bar W}_4}(6)$
maps $\fc (x)$ to one point, as we already know. 

Now suppose that $x$ is not regular. Since every abelian algebra
contaning $x$ is certainly included in the centralizer $\fc_0(x)$, 
we just need to understand the restriction of $\varphi$ to 
the linear subspace $\fc_0(x)$ to be able to determine the image 
of $x$ by $\varphi$.

If $x$ is not contained in ${\bar
X}_4\cup Y_4$, it has three eigenvalues, one of
which has multiplicity two. Up to conjugation, we may therefore suppose
that 
$$x=\begin{pmatrix} \a & 0 & 0 & 0 \\ 0 & \a & 0 & 0 \\
 0 & 0 & \b & 0 \\  0 & 0 & 0 & \g\end{pmatrix} \quad \in\quad
\fc (x)=\Bigg\{\begin{pmatrix} \d & \g & 0 & 0 \\ \eta & \varepsilon & 0 & 0 \\
 0 & 0 & \mu & 0 \\  0 & 0 & 0 & \nu\end{pmatrix}\Bigg\}.$$
Let $A$ denote the upper left corner of this matrix $M$ in $\fc (x)$. 
For $M$ to belong to (the cone over) ${\bar W}_4$, we have several
possibilities: either $\mu=\nu$, or $\mu$ or $\nu$ is an eigenvalue 
of $A$, or $A$ must be a homothety. This shows that the linear system
$I_{{\bar W}_4}(6)$, restricted to $\fc_0(x)$, contains a fixed hyperplane 
and two fixed quadrics, the residual system being generated by the three
linear conditions for $A$ to be a multiple of the identity. To 
resolve the indeterminacies we just need to blow-up the corresponding
codimension three linear subspace, and the image of 
$x$  by $\varphi$ is isomorphic with the fiber of that 
blow-up over $x$, which is just a projective plane. 

Now suppose that $x\in {\bar X}_4-X_4^0$. Then $x$ has two eigenvalues,
one of multiplicity three. The centralizer $c(x)$ is isomorphic with
$\fgl_3$, and the the linear system $I_{{\bar W}_4}(6)$ restricted to
$\fc(x)$ contains the fixed cubic $\det (M+tr(M)I)=0$. The residual
system is the space of cubics vanishing on the cone with vertex $I$
and base the variety of rank one matrices. This is the system of 
cubics on $\fsl_3$ vanishing on the projection of this rank variety, 
and we know by \cite{imcrelle} that its image is nothing but a copy
of the reduction variety $Red(3)$. \qed

\medskip\noindent {\it Remark}.
 This analysis suggests that it could be possible to
resolve the indeterminacies of $\varphi$ by blowing up successively
the different strata $X_4^0$, ${\bar X}_4$, $Y_4^0$, $Y_4$, $W_4$,
or rather their successive strict transforms -- but we have not 
been able to do that. 

Also it could be possible to extend this analysis to higher rank:
on each strata we can restrict the linear system that should define
$\varphi$, factor out the fixed components and get a linear system
that comes from smaller rank. This also makes sense for $a\ne 2$.

\subsection{The singular locus}

\begin{prop}
$Red(4)_{sing} = \cO_{min}\coprod \cO'_{min}$.
\end{prop}
 
\proof A simple computation shows that $\cO_6$ 
is contained in the regular locus (take local coordinates on the Grassmannian, 
write the commutativity conditions down and get $24$ independent 
linear relations). 
Since $\cO_6$  belongs to the closure of any orbit 
other than the two minimal orbits of dimension $3$, which we 
already know to be singular,
there is no other singular orbit. \qed

\medskip Recall that we denoted by $C_4=A_4$ the projectivized tangent
cone to a normal slice to $\cO_3'$ in $Red(4)$. This is an
eight-dimensional variety defined by $15$ quadratic equations. 
  
\begin{prop} 
The variety $C_4\subset\PP^{14}$ is projectively equivalent to $G(2,6)$. 
\end{prop}

\proof 
We define an equivariant  map $T$ from $\Lambda^2S^2U^*$ to the space 
of traceless symmetric 
maps from $U$ to $End(U)$, by sending an elementary  tensor $e^2\wedge
f^2$ 
to the map $B$ 
defined by 
$$B(u)(u)=(e,u)(f,u)\; e\we f, \quad u\in U,$$
with the identification of $\Lambda^2U^*$ with $U$. We claim that this
map 
$T$ sends 
the Grassmannian $G(2,S^2U^*)\subset\PP \Lambda^2S^2U^*$ isomorphically
on 
$C_4\subset 
\PP S_{1,0,-2}U$. 

Consider a generic point of $G(2,S^2U^*)$, that is, a generic pencil of 
conics in $\PP U\simeq\PP^2$. 
Such a pencil is defined by its base-locus, a set of four points in 
general position. 
Choosing homogeneous coordinates for which these four points are
$[1:0:0], 
 [0:1:0],  [0:0:1], [1:1:1]$,   
we get a pencil generated by the reducible conics $(x-z)y$ and
$(x-y)z$. 
But, by polarization, 
our map $T$ sends a tensor $ee'\wedge ff'$ to the map $B$ defined by 
$$ B(u)(u)=(e,u)(f,u)\; e'\we f'+(e',u)(f,u)\; e\we f'+(e,u)(f',u)\; 
e'\we f+(e',u)(f',u)\; e\we f.$$
Substituting $e=x-z, e'=y, f=x-y, f'=z$, we get 
$$B(u)=\begin{pmatrix} 3u_1 & -u_1 & -u_1 \\-u_2 & 3u_2 & -u_2 \\-u_3 & 
-u_3 & 3u_3 \end{pmatrix}
-(u_1+u_2+u_3)I.$$
Now for two vectors $u$ and $v$, let $\d_{ij}=u_iv_j-u_jv_i$. A simple 
computation shows that 
$$[B(u),B(v)]=\begin{pmatrix} \d_{12}+\d_{13} & -3\d_{12}+\d_{13} &
	       \d_{12}
-3\d_{13} \\
3\d_{12}+\d_{23} & -\d_{12}+\d_{23} & -\d_{12}-3\d_{23} \\
3\d_{13}-\d_{23} & 3\d_{23}-\d_{13} & -\d_{13}-\d_{23} 
\end{pmatrix},$$
and one can easily check that the image of this matrix is always
contained 
in $\langle u, v\rangle$. 
We conclude that $T$ maps $G(2,S^2U^*)$ to $C_4$, which are both 
irreducible of dimension $8$. 
Since $T$ is a linear automophism, it restricts to a projective 
equivalence between 
$G(2,S^2U^*)$ and $C_4$. \qed

\medskip By Lemma 10,$PGL_3$ has an open orbit in $G(2,6)=
G(2,S^2\CC^3)$, the space of pencils of plane conics. This is 
well-known and quite obvious, since a general pencil is determined 
by its base-locus -- four points in general position, and $PGL_3$
acts transitively on such four-tuples. In particular, we deduce
that the stabilizer of a general pencil is the stabilizer of 
its base-locus. This identifies for $n=4$ the finite group 
we introduced in Lemma 10:

\begin{coro} 
The finite group $K_4$ is the symmetric group $\cS_4$.
\end{coro}

\proof Given four points in general position in $\PP^2$, there is a 
unique projective 
transformation which fixes two of them and exchanges the other two. This 
implies that 
the stabilizer in $PGL_3$ of our four-tuple of points is a copy of
$\cS_4$, 
and the 
corresponding pencil of conics has the same stabilizer. \qed

\medskip But a more interesting consequence of the previous proposition
is:

\begin{coro} 
The variety of reductions $Red(4)$ is normal, with canonical 
singularities.
\end{coro}

\proof Since the Grassmannian $G(2,6)$ is projectively normal, the 
cone over it is normal, thus the tangent cone to a singular point 
of $Red(4)$ is normal as well. This implies that $Red(4)$ itself is
normal. By Theorem 8 its anticanonical divisor is $-K_{Red(4)}=
\cO_{Red(4)}(3)$, hence effective, and the singularities are
then automatically canonical. \qed  

\medskip\noindent {\it Remark}. As explained in \cite{katz}, $G(2,6)$
is also the projectivized tangent cone to a normal slice to the 
singular locus of $Hilb^4\PP^3$, which is also a $\PP^3$, 
parametrizing double points. What we expect is that the rational
map $\rho_4 : Hilb^4\PP^3\dashrightarrow Red(4)$ constructed in 
2.7, is a morphism contracting the divisor in $Hilb^4\PP^3$ defined
as the closure of linearly dependant four-tuples of points, to 
$\cO_3''\simeq {\check \PP}^3$, and restricting to an isomorphism 
outside this divisor, in particular around the singular locus, 
which should be mapped to $\cO_3'\simeq \PP^3$. Therefore the 
singularities should really be the same, and not just the tangent cones.

\subsection{Resolving the singularities}

Let ${\tilde G}$ denote the blow-up of $G(3,\fsl_4)$ along the 
smooth subvarieties $\cO_3'$ and $\cO_3''$. Since the tangent cone to
$Red(4)$ in a normal slice to each of these orbits is smooth, 
the strict transform of $Red(4)$ in ${\tilde G}$ is a smooth 
variety ${\tilde R}$ with an induced action of $PGL_4$. 
The two exceptional divisors are $G(2,6)$-fibrations above 
copies of $\PP^3$. 

Let $T$ denote a maximal torus in $PGL_4$. 

\begin{prop}
The smooth variety ${\tilde R}$ has only a finite number of fixed 
points of $T$. This number is equal to the Euler characteristic
$$\chi({\tilde R})=193.$$ 
\end{prop}

\proof A $T$-fixed point in  ${\tilde R}$ must dominate a $T$-fixed
point in $Red(4)$. Using our explicit description of the $PGL_4$-orbits
in $Red(4)$ we can easily determine these fixed points. Indeed, if we
choose for $T$ the torus defined by the canonical basis of $\CC^4$,
we see that an orbit $\cO$ contains a fixed point only when the 
corresponding representative is generated by diagonal matrices and 
matrices of the form $e_i^*\ot e_j$. Then all the fixed points in the
orbit can be deduced from a permutation of the basis vectors. 

We get the following numbers of fixed points in the different orbits:
$$\begin{array}{ccccccccccccccc}
\cO & \cO_{12} & \cO_{11} & \cO_{10}' & \cO_{10}'' & \cO_{9} & \cO_{8} 
& \cO_{8}' & \cO_{8}'' & \cO_{7} & \cO_{7}' & \cO_{7}'' & 
\cO_{6} & \cO_{3}' & \cO_{3}'' \\
\#\cO^T & 1 & 12 & 12 & 0 & 0 & 0 & 12 & 12 & 0 & 0 & 0 & 24 & 4 & 4
\end{array}$$
Each of these fixed points gives a unique fixed point in ${\tilde R}$,
except the eight ones in $\cO_{3}'\cup \cO_{3}''$. For each of these,
we need to count the number of normal directions that are fixed by
$T$ -- that is, the number of $T$ fixed points in the corresponding copy
of $G(2,6)$. It is easy to see that this number is finite, hence equal
to the Euler characteristic of the Grassmannian, that is $15$.
We thus get $120$ fixed points in ${\tilde R}$, plus $73$ coming 
from the smooth locus of $Red(4)$.

That the total number of fixed points equals the Euler characteristic
of ${\tilde R}$ is then an immediate consequence of the Byalinicki-Birula
decomposition \cite{bb}.
\qed

\begin{coro}
$Red(4)$ is rational. 
\end{coro}

\proof Since ${\tilde R}$ is smooth and has a finite number of points 
fixed by a torus action, it is a compactification of a $\CC^{12}$ - thus 
a rational variety, as well as $Red(4)$. \qed
 
\medskip The Byalinicki-Birula decomposition allows to compute the
Betti numbers of  ${\tilde R}$. For this we need the 
weights of the $T$-action on the tangent spaces to $\tilde{R}$ at 
the fixed points of $T$. 

For the $73$ fixed points that do not belong to the exceptional divisors
of the projection to $Red(4)$, we compute the tangent spaces to $\tilde{R}$
(or $Red(4)$, equivalently) as  limits 
of tangent spaces at points of the open $PGL_4$-orbit $\cO_{12}$. 
Indeed, the tangent space to $Red(4)$ at a point $\fa\in\cO_{12}$, as we 
have seen, is easily computed as the image of the (injective) map 
$$\fsl_4/\fa\ra Hom(\fa,\fsl_4/\fa)=T_{\fa}G(3,\fsl_4)$$ defined by 
the Lie bracket. Note that we need only one computation per 
$PGL_4$-orbit, 
since the symmetric group $\cS_4$ acts transitively on the set of
$T$-fixed points in each orbit. Thus only six computations are enough
to take care of these $73$ fixed points.  
 
For the $120$ remaining fixed points, we proceed as follows.
Consider the point $\fa$ of $\cO_3''$ defined as at the beginning of 2.5,
with $n=4$. The splitting of $\CC^4$ into the sum of the hyperplane
$U$ and the line $\ell$ generated by $e_4$ leads to the identifications 
$$\begin{array}{ccc}
T_{\fa}G(3,\fsl_4) & \simeq & Hom(\ell^*\ot U, U^*\ot U\op U^*\ot\ell) \\
\cup & & \cup \\
T_{\fa}Red(4) & \simeq & Hom^s(\ell^*\ot U, U^*\ot U)
\end{array}$$
where $Hom^s(\ell^*\ot U, U^*\ot U):=\ell\ot S^2U^*\ot U\subset
\ell\ot U^*\ot U^*\ot U=Hom(\ell^*\ot U, U^*\ot U)$. Now, 
recall that $S^2U^*\ot U = U^*\op S_{1,0,-2}U$. The $U^*$ factor
corresponds to the tangent directions to the orbit $\cO_3''$. 
The other term $S_{1,0,-2}U=\we^2(S^2U^*)\ot\det U$ is, up to  a
twist, the ambient space for the Pl\"ucker embedding of $G(2,S^2U^*)$,
which we identified with the projectivized 
tangent cone to $Red(4)$ in the directions
normal to $\cO_3''$. Then the fixed points of $T$ in $\tilde{R}$
over this point $\fa$ of $Red(4)$, are in correspondence with the 
$15$ fixed points of $T$ contained in that Grassmannian. And we 
deduce the weights of the $T$-action on the tangent space to 
$\tilde{R}$ from those of the $T$-action on the tangent space
to $G(2,S^2U^*)$, through the previous identifications.  
Again, there are enough symmetries for the effective 
computations to remain tractable. 

Finally, we choose a general enough one-dimensional subtorus 
of $T$, and count the number of negative weights of the restricted
action on the tangent spaces to the fixed points : this
gives the dimensions of the corresponding strata in the 
Byalinicki-Birula decomposition. 
The conclusion is the following: 

\begin{prop}
The odd Betti numbers of  ${\tilde R}$ are all zero. The
even Betti numbers are
$$1,3,9,15,23,29,33,29,23,15,9,3,1.$$
\end{prop}

Applying the same arguments as for the proof of Theorem 2.4 in 
\cite{katz}, we can deduce the ranks of the Chow groups of $Red(4)$.
Indeed, passing from $Red(4)$ to $\tilde{R}$ amounts to replacing
two copies of $\PP^3$ by two $G(2,6)$-bundles over them, and the 
ranks of the Chow groups are modified accordingly. We get:
 
\begin{prop}
The Chow groups of $Red(4)$ have respective ranks 
$$1,1,3,5,7,11,14,13,11,7,5,1,1.$$
In particular, $Red(4)$ has Picard number one.
\end{prop}

\bigskip\noindent 
Atanas ILIEV,
 
\noindent 
Institute of Mathematics, Bulgarian Academy of Sciences, 

\noindent 
Acad. G. Bonchev street 8, 1113 Sofia, Bulgaria.

\noindent Email : ailiev@math.bas.bg

\smallskip\noindent 
Laurent MANIVEL, 

\noindent 
Institut Fourier, UMR 5582 (UJF-CNRS), 

\noindent 
BP 74, 38402 St Martin d'H\`eres Cedex, France.

\noindent Email : Laurent.Manivel@ujf-grenoble.fr

\end{document}